\documentclass[12pt,twoside,reqno]{amsart}
\usepackage{amsxtra,amscd}
\usepackage{graphicx}
\usepackage{amsmath}
\usepackage{amsfonts}
\usepackage{amssymb}

\newtheorem{theorem}{Theorem}[section]

\theoremstyle{definition}
\newtheorem{definition}[theorem]{Definition}

\newtheorem{remark}[theorem]{Remark}

\numberwithin{equation}{section}

\begin{document}
\title{Notes on the section conjecture of Grothendieck}
\author{Feng-Wen An}
\address{School of Mathematics and Statistics, Wuhan University, Wuhan,
Hubei 430072, People's Republic of China}
\email{fwan@amss.ac.cn}

\begin{abstract}
In this short note, we will give the key point of the section conjecture of Grothendieck, that is reformulated by monodromy actions. Here, we will also give the result of the section conjecture for algebraic schemes over a number field.
\end{abstract}

\maketitle

\begin{center}
{\tiny {Contents} }
\end{center}

{\tiny \qquad {Introduction} }

{\tiny \qquad {1. Preliminaries} }

{\tiny \qquad {2. Section Conjecture: Key Point} }

{\tiny \qquad {3. Section Conjecture, II: Algebraic Schemes}}

{\tiny \qquad {References}}

\section{Preliminaries}

\subsection{Notation}

In this note, an \textbf{algebraic $K-$variety} is an integral scheme over a number field $K$ of finite type.

An \textbf{arithmetic variety} is an integral scheme $X$
satisfying the conditions:

\begin{itemize}
\item $\dim X \geqq 1$.

\item There is a surjective morphism $f:X \to Spec\left( \mathbb{Z}\right) $
of finite type.
\end{itemize}

For an integral scheme $Z$, set

\begin{itemize}
\item $k(Z)\triangleq$ the function field of an integral scheme $Z;$

\item $\pi _{1}^{et}\left( Z\right) \triangleq$ the \'{e}tale fundamental
group of $Z$ for a geometric point of $Z$ over a separable closure of the
function field $k\left( Z\right).$
\end{itemize}

In particular, for a field $L$, we set $$\pi _{1}^{et}(L)\triangleq\pi _{1}^{et}(Spec(L)).$$

\subsection{Outer homomorphisms}

Let $G,H,\pi_{1},\pi_{2}$ be four groups with homomorphisms $p:G\to \pi_{1}$ and $q:H\to \pi_{2}$, respectively.

The \textbf{outer homomorphism set} from $G$ into $H$ over $\pi_{1}$ and $\pi_{2}$ respectively, denoted by $Hom_{\pi_{1},\pi_{2}}^{out}(G,H)$, is the set of the maps $\sigma$ from the quotient $\frac{\pi_{1}}{p(G)}$ into the quotient $\frac{\pi_{2}}{q(H)}$ given by a group homomorphism $f:G\to H$ in such a manner: $$\sigma:\frac{\pi_{1}}{p(G)} \to \frac{\pi_{2}}{q(H)}, x\cdot p(G)\mapsto f(x)\cdot q(H)$$ for any $x\in \pi_{1}$.

\begin{remark}
Let $Out(G)$ and $Out(H) $ be the outer automorphism groups. Suppose that $G$ and $H$ are normal subgroups of $\pi_{1}$ and $\pi_{2}$, respectively.
Then $Hom_{\pi_{1},\pi_{2}}^{out}(G,H)$ can be regarded as a subset of $Hom(Out(G),Out(H))$.
\end{remark}

\subsection{Universal Cover and \emph{sp}-Completion}

For details, see \cite{An3}-\cite{An8}.

\begin{definition}
Let $X$ be an arithmetic variety (or  algebraic $K-$variety, respectively). Then there is an integral variety $X_{\Omega
_{et}}$ and a surjective morphism $p_{X}:X_{\Omega _{et}}\rightarrow X$
satisfying the conditions:

\begin{itemize}
\item $k\left( X_{\Omega _{et}}\right) ={k(X)}^{un}$ (or $={k(X)}^{au}$, respectively);

\item $p_{X}$ is affine;

\item $k\left( X_{\Omega _{et}}\right) $ is Galois over $k\left( X\right) ;$

\item $X_{\Omega _{et}}$ is quasi-galois closed over $X$ by $p_{X}$.
\end{itemize}

The integral variety $X_{\Omega _{et}}$ is called a \textbf{universal
cover} over $X$ for the \'{e}tale fundamental group $\pi _{1}^{et}\left( X\right) $, denoted by $
\left( X_{\Omega _{et}},p_{X}\right) .$
\end{definition}

\begin{definition}
For any integral variety $X$, there exists an integral variety $X_{sp}$ and a surjective morphism $\lambda_{X}:X_{sp}\to X$ such that
\begin{itemize}
\item $\lambda_{X}$ is affine;

\item $X_{sp}$ is $sp-$complete;

\item $k(X_{sp})$ is an algebraic closure of $k(X)$;

\item $X_{sp}$ is quasi-galois closed over $X$ by $\lambda_{X}$.
\end{itemize}

 The integral scheme $X_{sp}$, is said to be an \textbf{$sp-$completion} of $X$, denoted by $(X_{sp},\lambda_{X})$.
\end{definition}

\subsection{Formally unramified extensions}

(See \cite{An6,An9}).
Let us recall the definition for formally unramified extension of a given field.

\begin{definition}
Let $K_{1}$ and $K_{2}$ be two arbitrary extensions over a field $K$ such
that $K_{1}\subseteq K_{2}$.

$\left( i\right) $ $K_{2}$ is said to be a \textbf{finite formally unramified Galois}
extension of $K_{1}$ if there are two algebraic varieties $X_{1}$ and $X_{2}$
and a surjective morphism $f:X_{2}\rightarrow X_{1}$ such that

\begin{itemize}
\item $k\left( X_{1}\right) =K_{1},k\left( X_{2}\right) =K_{2}$;

\item $X_{2}$ is a finite \'{e}tale Galois cover of $X_{1}$ by $f$.
\end{itemize}

$\left( ii\right) $ $K_{2}$ is said to be a \textbf{finite formally unramified}
extension of $K_{1}$ if there is a field $K_{3}$ over $K$ such that $K_{2}$
is contained in $K_{3}$ and $K_{3}$ is a finite formally unramified Galois extension
of $K_{1}.$

$\left( iii\right) $ $K_{2}$ is said to be a \textbf{formally unramified} extension
of $K_{1}$ if the field $K_{1}(\omega)$ is a finite formally unramified extension of $
K_{1}$ for each element $\omega \in K_{2}$. In such a case, the element $
\omega$ is said to be \textbf{formally unramified} over $K_{1}$.
\end{definition}

Let $L$ be an arbitrary extension over a field $K$. Set

\begin{itemize}
\item $L^{al}\triangleq $ an algebraical closure of $L$;

\item $L^{au}\triangleq $ the mximal finite formally unramified
subextensions over $L$ contained in $L^{al}$;

\item $G(L)\triangleq$ the absolute Galois group $Gal(L^{al}/L)$;

\item $G(L)^{au}\triangleq$ the Galois group $Gal(L^{au}/L)$ of the maximal formally unramified extension $L^{au}$ of $L$.
\end{itemize}

\begin{remark}
It is seen that for the case of an algebraic extension, in general, the formally unramified extension defined in \emph{Definition 1.4} does not coincide  with that in algebraic number theory.
\end{remark}

\begin{remark}
Note that we define another unramified extensions in \cite{An4*,An5,An8} for arithmetic schemes, which is a generalization of unramified extensions in algebraic number theory and hence is different from the above one defined in \emph{Definition 1.4}.
\end{remark}

\section{Section Conjecture: Key Point}

The key point of the section conjecture, or what Grothendieck tells to us indeed, can be reformulated in a manner such as the following.

\begin{theorem}
Let $X,Y$ be two algebraic $K-$varieties or two arithmetic varieties.
Then there exists a bijection $$Hom(X,Y) \cong Hom(\frac{Aut(X_{sp}/X)}{Aut(X_{\Omega _{et}}/X)},\frac{Aut(Y_{sp}/Y)}{Aut(Y_{\Omega _{et}}/Y)})$$ between sets.
\end{theorem}

This result is immediate from the monodromy actions of automorphism groups on the integral schemes that are considered. See \cite{An8,An9} for the proof.

Here, we have the following group isomorphisms $$Gal(k(X_{sp})/k(X))\cong Aut(X_{sp}/X);$$ $$Gal(k(Y_{sp})/k(Y))\cong Aut(Y_{sp}/Y);$$ $$\pi_{1}^{et}(X)\cong Gal(k(X_{\Omega _{et}})/k(X))\cong Aut(X_{\Omega _{et}}/X);$$  $$\pi_{1}^{et}(Y)\cong Gal(k(Y_{\Omega _{et}})/k(Y))\cong Aut(Y_{\Omega _{et}}/Y)$$ in virtue of the properties of quasi-galois closed schemes.

\section{Section Conjecture, II: Algebraic Schemes}

For anabelian geometry of algebraic schemes over $K$, we have the following results. See \cite{An9} for detail.

\begin{theorem}
Let $X$ and $Y$ be two algebraic $K-$varieties such that $k\left( Y\right) $ is contained in $ k\left( X\right)$. Then there is a bijection
\begin{equation*}
Hom\left( X,Y\right) \cong Hom_{\pi _{1}^{et}(k(X)),\pi _{1}^{et}(k(Y))
}^{out}\left( \pi _{1}^{et}\left( X\right) ,\pi _{1}^{et}\left( Y\right)
\right)
\end{equation*}
between sets.
\end{theorem}

\begin{theorem}
For any algebraic $K-$variety $X$, there
is a bijection
\begin{equation*}
\Gamma \left( X/K\right) \cong Hom_{\pi _{1}^{et}(K),\pi _{1}^{et}(k(X)) }^{out}\left( \pi _{1}^{et}(K) ,\pi _{1}^{et}\left( X\right) \right)
\end{equation*}
between sets.
\end{theorem}

Using Galois groups of fields,  we will have the following versions of the main theorems above.

\begin{theorem}
Let $X$ and $Y$ be two algebraic $K-$varieties such that $k\left( Y\right) $ is contained in $ k\left( X\right)$. Then there is a bijection
$$
Hom\left( X,Y\right)  \cong Hom_{G(k(X)),G(k(Y))
}^{out}\left( G(k(X))^{au} ,G(k(Y))^{au}
\right)
$$
between sets.
\end{theorem}

\begin{theorem}
For any algebraic $K-$variety $X$, there
is a bijection
\begin{equation*}
\Gamma \left( X/K\right) \cong Hom_{G(K),G(k(X)) }^{out}\left( G(K)^{au} ,G\left( k(X)\right)^{au} \right)
\end{equation*}
between sets.
\end{theorem}


\newpage

\begin{thebibliography}{99}

\bibitem{An1} An, F-W. The combinatorial graph of a scheme. eprint arXiv:0801.2609.

\bibitem{An2} An, F-W. Automorphism groups of quasi-galois closed arithmetic
schemes. eprint arXiv:0907.0842.

\bibitem{An2*} An, F-W. Automorphism groups of quasi-galois closed
arithmetic schemes, (without affine structures). preprint.

\bibitem{An3} An, F-W. On the existence of geometric models for function
fields in several variables. eprint arXiv:0909.1993.

\bibitem{An4*} An, F-W. on the \'{e}tale fundamental groups of arithmetic
schemes, revised. eprint arXiv:0910.4646.

\bibitem{An5} An, F-W. On the arithmetic fundamental groups. eprint
arXiv:0910.0605.

\bibitem{An6} An, F-W. On the algebraic fundamental groups. eprint
arXiv:0910.4691.

\bibitem{An7} An, F-W. Notes on the quasi-galois closed schemes. eprint
arXiv:0911.1073.

\bibitem{An8} An, F-W. On the section conjecture of Grothendieck.
eprint arXiv:0911.1523.

\bibitem{An9} An, F-W. On the section conjecture of Grothendieck, II.
preprint.

\bibitem{F-K} Freitag, E; Kiehl, R. \'{E}tale Cohomology and the Weil
Conjecture. Springer, Berlin, 1988.

\bibitem{EGA} Grothendieck, A; Dieudonn\'{e}, J. \'{E}l\'{e}ments de G\'{e}%
oem\'{e}trie Alg\'{e}brique. vols I-IV, Pub. Math. de l'IHES, 1960-1967.

\bibitem{SGA1} Grothendieck, A; Raynaud, M. Rev$\hat{e}$tements $\acute{E}$%
tales et Groupe Fondamental (SGA1). Springer, New York, 1971.

\bibitem{faltings} Grothendieck, A. Letter to Faltings, in \emph{Geometric
Galois Actions}, Vol 1, edited by Schneps, L ; Lochak, P. Cambridge
University Press, New York, 1997.

\bibitem{Hrtsh} Hartshorne, R. Algebraic Geometry. Springer, New York, 1977.

\bibitem{Mln} Milne, J. \'{E}tale Cohomology. Princeton University Press,
Princeton, New Jersey, 1980.

\bibitem{pop} Pop, F. Glimpses of Grothendieck's anabelian geometry, in
\emph{Geometric Galois Actions}, Vol 1, edited by Schneps, L ; Lochak, P.
Cambridge University Press, New York, 1997.

\bibitem{sv1} Suslin, A; Voevodsky, V. Singular homology of abstract
algebraic varieties. Invent. Math. 123 (1996), 61-94.

\bibitem{sv2} Suslin, A; Voevodsky, V. Relative cycles and
Chow sheaves, in \emph{Cycles, Transfers, and Motivic Homology
Theories}, Voevodsky, V; Suslin, A; Friedlander, E M. Annals of Math
Studies, Vol 143. Princeton University Press, Princeton, NJ, 2000.

\bibitem{www} http://people.math.jussieu.fr/~leila/grothendieckcircle/letters.php
\end{thebibliography}
\end{document}